# Hierarchical Electricity and Carbon Trading in Transmission and Distribution Networks Based on Virtual Federated Prosumer

Lu Wang, *Student Member, IEEE*, Zhi Wu, *Member, IEEE*, Wei Gu, *Senior Member, IEEE*, Haifeng Qiu, *Student Member, IEEE*, Shuai Lu, *Graduate Student Member, IEEE*

*Abstract*—**Facing the dilemma of growing energy demand and mitigating carbon emissions, this paper proposes an energy sharing mechanism based on virtual federated prosumers (VFPs) with budget allocation for joint electricity and carbon market to incentivize distributed energy resources to participate in the hierarchical market and reduce carbon emissions. At the transmission level, the regional transmission operator coordinates transactions between two markets, the inter-VFP energy sharing market and the wholesale market, intending to minimize the overall cost of VFPs. The energy sharing market clearing problem is formulated as a generalized Nash game, for which we develop a first-order response algorithm to obtain the equilibrium. At the distribution level, the VFPs play the role of selfless auctioneer that leverage discriminatory weights and benchmark prices to allocate the electricity-carbon budget among entities in the VFP to maximize social welfare. The Nash game is exploited to characterize the budget allocation problem, for which a distributed feedback allocation algorithm is proposed. The entire hierarchical electricity and carbon trading is modeled as an equilibrium problem and is solved iteratively. Case studies based on a practical regional grid verify the effectiveness of the proposed algorithm and show that the mechanism is effective in improving energy efficiency and reducing carbon emissions.**

*Index Terms*--Carbon emission, distributed energy resources, virtual federated prosumer, energy market, budget allocation.

## NOMENCLATURE

### A. Abbreviation

| | |
|---|---|
| VPP | Virtual power plant |
| VFP | Virtual federated prosumers |
| GNG | Generalized Nash game |
| RTO | Regional transmission operator |
| DER | Distributed energy resources |
| DC | Direct current |
| DG | Distributed Generation |
| ESS | Energy storage stations |
| PV | Photovoltaic station |
| GT | Gas turbine |
| NE | Nash equilibrium |
| FRA | First-order response algorithm |
| DFAA | Distributed feedback allocation algorithm |

### B. Sets and Indices

| | |
|---|---|
| $\mathcal{I}$ | Set of indices of VFPs, indexed by $i$ |
| $\mathcal{N}$ | Set of indices of consumers, indexed by $n$ |
| $\mathcal{M}$ | Set of indices of producers, indexed by $n$ |
| $L$ | Set of indices of power lines |
| $\widehat{E_n}$ | Set of physical constraints of consumers |
| $\widehat{E_m}$ | Set of physical constraints of producers |
| $\mathcal{P}_i$ | The strategy set of VFP, indexed by $i$ |
| $\mathcal{P}_{I+1}$ | The strategy set of RTO, indexed by $I+1$ |
| $\mathcal{R}$ | Payoff function at the game $\mathcal{G}$ |
| $\mathcal{F}$ | Payoff function at the game $\Theta(\mathcal{N}, \Omega_N, \mathcal{F})$ |
| $Z$ | Payoff function at the game $\Theta(\mathcal{M}, \Omega_M, Z)$ |
| $h_i$ | Cost function of VFP $i$ production electricity |
| $\Omega_n / \Omega_m$ | The strategy set of consumers/producers |

### C. Variables

| | |
|---|---|
| $t$ | Index of periods, from 1 to 24 |
| $\lambda_s, \lambda_c$ | Price of shared electricity/carbon at transmission level |
| $\lambda_s^*, \lambda_c^*$ | Price of shared electricity/carbon at the equilibrium point |
| $P_{i,s/i,c}$ | Share power/carbon emissions from VFP $i$ |
| $P_{i,w}$ | Electricity purchased by VFP $i$ from the wholesale market |
| $P_{i,g}$ | Internal electricity production of VFP $i$ |
| $P_{s/c}^*$ | Share power/carbon emissions at the equilibrium point |
| $P_w^*$ | Electricity purchased from the wholesale market at the equilibrium point |
| $P_{i,j}^t$ | Transmission power of line between buses $i$ and $j$ at $t$ |
| $P_{i,DG/GT}^t$ | Output electrical power of DG/GT in VFP $i$ at $t$ |
| $P_{i,dis/cha}^t$ | The charging and discharging power of ESS in VFP $i$ at $t$ |
| $P_{i,bat}^t$ | State of charge of ESS in VFP $i$ at $t$ |
| $B_{i,j}^t$ | Conductance between lines $i$ and $j$ |
| $\theta_{ij}^t$ | Voltage phase angle of node $i$ |
| $p_i, q_j$ | Price weight of power sent by VFP to consumer $n$ / producer $m$ |
| $e_{n/m/\psi}$ | Energy demand / supply / surplus of consumer $n$ / producer $m$ / prosumer $\psi$ |
| $\kappa_{n/m/\psi}$ | Energy price of consumer $n$ / producer $m$ / prosumer $\psi$ |
| $\Gamma_i$ | Budget for VFP $i$ |
| $\mu_n, \mu_m$ | Benchmark clearing price per unit of consumer/producer by VFP (¥/MWh) |
| $D_{i,in}, D_{i,out}$ | The 0-1 variables indicate the state of VFP interaction with the upstream grid |
| $u_i^t$ | The 0-1 variables indicate the state of GT |

| $\Upsilon^*$ | Budget allocation at NE of the game $\Theta(\mathcal{N}, \Omega_N, \mathcal{F})$ |
| --- | --- |
| $z^*$ | Budget allocation at NE of game $\Theta$ |
| $c_{co2/m}$ | Cost of carbon/producer $m$ |

### D. Parameters

| $P_{i,j}^{max/min}$ | Upper and lower limits of active capacity between transmission lines $i$ and $j$ |
| --- | --- |
| $P_{i,GT}^{max/min}$ | Upper and lower limits of GT output power |
| $P_{i,DG}^{max/min}$ | Upper and lower limits of DG output power |
| $P_{i,bat}^{max/min}$ | Upper and lower capacity limits in ESS |
| $\eta_{i,cha/dis}$ | Charging and discharging efficiency of ESS |
| $P_{i,w}^{max/min}$ | Upper and lower limits of power purchased from the wholesale market |
| $P_c^{max}$ | Maximum intra-day carbon emissions in the region |
| $\varepsilon_i$ | Carbon emission factors of VFP $i$ |
| $\alpha_i, \beta_i$ | Electricity inverse demand function coefficient of VFP $i$ |
| $\omega_n, \upsilon_n$ | Quadratic and primary term coefficients of the utility function of consumer $n$ |
| $a_m, b_m$ | Quadratic and primary term coefficients of the cost function of producer $m$ |
| $M, N, I$ | Total number of producers/consumers/VFPs |
| $\delta, \Delta l$ | Convergence step |
| $\epsilon$ | Convergence tolerance |

## I. INTRODUCTION

### A. Background and Motivation

THE sustainable development of human society is restricted by climate change and energy shortage [1]. Many countries around the world are vigorously cascading the utilization of renewable energy while incenting emission reduction. With the prevalence of distributed energy resources (DER), the consumers passively participating in the power market have gradually transformed into proactive prosumers [2], and various innovative trading methods are proposed for the participation of prosumers in the power market.

The market trade mechanisms for the prosumers usually change with the scale of the prosumer. Nevertheless, small-scale prosumers cannot directly trade energy at fair prices in the wholesale market due to their small capacity. Instead, they share surplus energy via a peer-to-peer (P2P) platform or sign a contract with retailers to participate in the wholesale market [3]. Although the P2P platform provides a market environment for prosumers to participate in power trading, the efficient allocation and utilization of DER capacity cannot be guaranteed [4]. In addition, another way for prosumers to participate in the market is to aggregate a large number of DER to conduct market transactions. In such cases, DERs are treated as tradable virtual resources, such as virtual microgrids, load aggregators, and virtual power plants (VPP). However, the consumers aggregated by VPP are treated as price-taker without discrimination due to the top-down structure of VPP, which ignores rational individual preferences [4]. Moreover, the incentive results of VPP are not consistent with the strategies based on environmental considerations. Hence, it is crucial to combine the electricity and carbon markets, considering the mutual impact of the market and the environment [5], to effectively restructure energy demand and promote a low-carbon transition.

The future low-carbon energy market, where consumers, producers and consumers can all participate at a fair and reasonable price, requires a new market mechanism and pricing scheme in the transmission and distribution network.

### B. Literature Review

Recently, energy sharing mechanisms have been proposed and widely utilized in power systems for compatibility with the increasing number of small-scale DERs within the system. Such small-scale DERs are usually involved in energy sharing at the distribution network level as prosumers which is divided into two main categories depending on their role in the market, the first being as price-takers. Ref. [6] characterized the interaction between prosumers and grid company in terms of a one-leader, N-follower Stackelberg game for realizing P2P trading market clearing. Furthermore, Paudel et al. [7] employed a more complex N-leader, M-follower Stackelberg game model for market clearing. More realistically, [8] considered user preferences in the modeling to encourage the use of clean energy. Moreover, the privacy of the gamers was protected while enabling energy sharing in [9]. For the second category, the prosumer participates in the market as a price-taker. A bilateral auction mechanism was presented in [10] and an adaptive pricing strategy was designed. Chen. et al. adopted a series of models developed by the GNG to characterize energy sharing markets with a large number of prosumers in [11]-[12]. Besides these two approaches, it is also a common approach to utilize allocation mechanisms in shared energy markets. Ref. [13] and Ref. [14] modeled the allocation of costs and computational burdens respectively with the use of non-cooperative games. In contrast, cooperative games are utilized to characterize the process of energy sharing in [15]-[16].

However, most of the above mechanisms are restricted to the independent microgrid or distribution network level, where DERs cannot participate directly in the upper wholesale market except through virtual federations. The authors of [17] developed a virtual microgrid through the reorganization of communication systems to improve the transaction efficiency and social welfare of prosumers in the region. Similarly, [18] proposed that virtual energy hubs participate in multiple types of markets to minimize their costs. The VPP is proposed to cooperate DER effectively, allowing DER to get indirectly involved in the wholesale market without specific geographical location and composition restrictions. The bidding strategies of VPPs in the energy market [19], DER dispatching [20] and system uncertainty have been well studied [21].

Although DERs can participate in different levels of market transactions through these virtual federations, none of the above studies consider the impact of market mechanisms on carbon emissions, which is incompatible with the urgent goal of reducing carbon emissions today. Most of the existing research on electricity and carbon markets has focused on planning [22], dispatching [23], and optimizing [24] the grid at a single level with the goal of low carbon emissions. The pricing of carbon emissions consists of two approaches: carbon tax [25] and cap-and-trade system [26]. Pricing with a carbon tax is vulnerable to energy price volatility. Comparatively, the latter is a more intuitive reflection of the level of control of carbon emission reduction which is particularly suitable for implementation in the small-scale region [27]. Nevertheless, the above studies all focus on carbon trading markets established for a single tier of the grid, which cannot be adapted to a hierarchical market with extensive DER access.

As a result, to build such an ideal energy market, three challenges exist, they are (i) how DER's participation in wholesale markets can improve the value of energy use while getting a fair price; (ii) designing incentive-compatible allocation mechanisms for individually rational and flexible consumers; (iii) how to clear the joint power and carbon markets to reduce carbon emissions in a transmission and distribution hierarchical market structure.

### C. Contribution and Organization

Given this context, this paper proposes an energy sharing mechanism with budget allocation in the joint electricity and carbon market to incentivize DERs and flexible consumers to participate in the market and reduce carbon emissions. The contributions of this paper are as follows:

1) A novel market stakeholder of "virtual federated prosumers (VFPs)", which spontaneously formed by flexible consumers, traditional producers, and DER, is proposed to incentive players to participate in shared energy markets and wholesale markets. Meanwhile, a VFP-based hierarchical electricity and carbon trading framework is established, which includes an energy sharing mechanism at the transmission level and budget allocation at the distribution level to reduce carbon emissions.

2) An energy sharing mechanism is devised to coordinate the energy among VFP and wholesale market at the transmission level, which is then captured by a generalized Nash game (GNG) model, and a first-order response algorithm is proposed to obtain the equilibrium. In contrast to the existing energy sharing mechanism which is mostly limited to within independent microgrids, we extend it to regional grids based on existing research, prioritizing intra-regional energy trading while connecting to wholesale markets, thus effectively promoting intra-regional DER consumption.

3) An electricity-carbon budget allocation approach based on the Nash game is proposed to integrate the selfish interests of the rational individual with the overall profit of VFP at the distribution level, which improves social welfare. A distributed feedback allocation algorithm is proposed for budget allocation to reduce costs and intra-regional carbon emissions.

The organization of the remainder of this article is as follows: In Section II, an energy sharing mechanism with a budget allocation is presented. In Section III, the game model is established, and the existence and uniqueness of the equilibrium are proved. A distributed algorithm is proposed in Section IV and case studies are presented in Section V. In Section VI, the paper is concluded.

## II. ENERGY SHARING MECHANISMS WITH BUDGET ALLOCATION

### A. Assumptions

A hierarchical electricity and carbon sharing mechanism between transmission and distribution networks is focused on in this paper. The assumptions are as follows:

a) VFPs are virtual players formed spontaneously by consumers, traditional polluting producers and green producers (containing PVs or ESSs). At the transmission level, VPF plays the role of bargainer for energy sharing and participates in the wholesale market, whereas at the distribution level, VFP acts as the selfless and trustful center to allocate the budget.

b) The Regional Transmission Operator (RTO) serves as a non-profit synergy center between the transmission network and the upstream grid to clear the upper energy sharing market.

c) The government agency will announce the next day's 24-hour carbon credits within the different VFP areas the next day.

d) The cost function of VFP $i$ $h_i$ are differentiable and convex. All the supplied power $P_{i,g} + P_{i,w}$ in the system at each moment meets the demand $\sum_{i \in} P_{i,load}$.

### B. Problem Description

The framework of joint electricity and carbon markets is divided into transmission-level and distribution-level, as shown in Fig. 1.

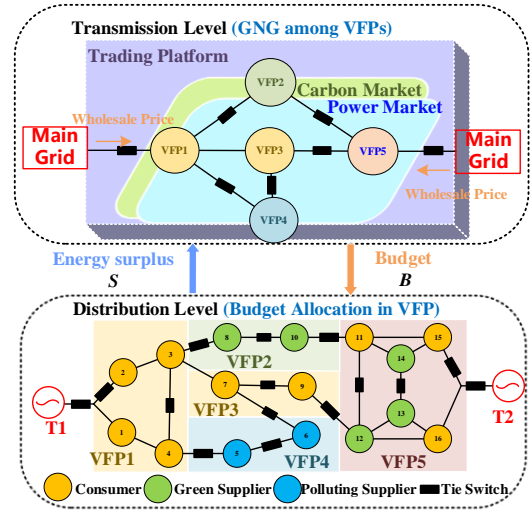

Fig. 1. The framework of the power and carbon market at the transmission and distribution level

At the transmission level, shared energy market players consist of VFPs, RTO and upstream power grid. The RTO acts as a particular player, coordinating stakeholders to determine shared electricity prices $\lambda_s$ and carbon prices $\lambda_c$ based on a shared electricity $P_{i,s}$ and carbon emissions $P_{i,c}$ among VFPs. VFP $i \in \mathcal{I} = \{1, ..., I\}$ determines a trading strategy to minimize costs based on the shared electricity price $\lambda_s$, carbon price $\lambda_c$ and wholesale market price $\lambda_w$, and will finalize the amount of shared electricity $P_{i,s}$, carbon emissions $P_{i,c}$ and trading with the wholesale market $P_{i,w}$. After multiple iterations of interaction between VFPs and RTOs to reach equilibrium, the shared energy market clearing at the transmission network layer is completed, at which point every $VFP_i$ gets its electricity-carbon budget $\Gamma_i$, which can be cost, revenue, or profit, depending on the specific type of VFP. The whole process of clearing the energy sharing market can be characterized by a generalized Nash game (GNG), due to the common constraints among the players, i.e., the balance of supply and demand.

In terms of the distribution level, the consumers indexed by $n \in \mathcal{N}_i := \{1, \ldots, N\}$ and distributed energy producers $m \in \mathcal{M}_i := \{1, \ldots, M\}$ spontaneously formed the $VFP_i$, which can be classified into different types, including federated consumers, federated producers, federated green prosumers and federated polluted prosumers. These small-scale prosumers with DER participate in the regional energy market via VFP to maximize their profits, which will lead to competition for budget shares. In this case, VFP acts as the center to allocate the budget $\Gamma_i$ within the coalition to achieve the internal social welfare optimum. Nodes within the VFP $n/m$ compete on a budget based on their utility to determine their energy demand/supply $e_{n/m}$. In this process of budget competition, the role of the players can be ignored and all can be captured by the Nash game.

Then, according to the $E_i^{k+1} := \sum_{\mathcal{N}_i / \mathcal{M}_i} n/m$, VFPs are engaged in energy sharing with each other again at the transmission level to strive for higher profits until they get a stable budget. Eventually, the market equilibrium is achieved through multiple iterations of transmission and distribution networks.

### C. VFP Sharing Energy Model at Transmission Level

At the transmission-level, the objective function of VFP $i \in \mathcal{I} = \{1, \ldots, I\}$ is given in (1a), which contains costs of electricity generation, costs of shared electricity between VFPs, costs of interaction with wholesale markets, and costs of carbon emissions.

$$\min_{P_{i,s}, P_{i,c}, P_{i,w}} \left( \begin{array}{c} h_i(P_{i,g}) + \sum_{t \in T} \lambda_s^t \cdot P_{i,s}^t + \sum_{t \in T} P_{i,w}^t \cdot D_{i,in}^t \cdot \lambda_{i,w,in}^t \\ + \sum_{t \in T} \lambda_c^t \cdot P_{i,c}^t - \sum_{t \in T} P_{i,w}^t \cdot D_{i,out}^t \cdot \lambda_{i,w,out}^t \end{array} \right) \quad (1a)$$

$$\text{s.t.} \sum_{(i,j) \in L(i)} P_{i,j}^t = \sum_{i \in \mathcal{I}} P_{i,GT}^t + \sum_{i \in \mathcal{I}} P_{i,DG}^t - P_{i,Load}^t$$
$$+ \sum_{i \in \mathcal{I}} P_{i,w}^t + \sum_{i \in \mathcal{I}} (P_{i,dis}^t - P_{i,cha}^t) \quad (1b)$$

$$B_i^t (\theta_i^t - \theta_j^t) = P_{i,j}^t \quad (1c)$$

$$\theta_{slack} = 0 \quad (1d)$$

$$P_{i,j}^{t,\min} \leq P_{i,j}^t \leq P_{i,j}^{t,\max}, \quad \forall i,j \in \mathcal{I}, i \neq j \quad (1e)$$

$$u_i^t P_{i,GT}^{t,\min} \leq P_{i,GT}^t \leq u_i^t P_{i,GT}^{t,\max}, \forall i \in \mathcal{I} \quad (1f)$$

$$P_{i,DG}^{t,\min} \leq P_{i,DG}^t \leq P_{i,DG}^{t,\max}, \quad \forall i \in \mathcal{I} \quad (1g)$$

$$\begin{cases} P_{i,bat}^{\min} \leq P_{i,bat}^t \leq P_{i,bat}^{\max} \\ P_{i,bat}^{\min} \leq P_{i,bat}^{t-1} \leq P_{i,bat}^{\max} \\ P_{i,bat}^t = P_{i,bat}^{t-1} + P_{i,cha}^t \eta_{i,cha} - P_{i,dis}^t / \eta_{i,dis} \end{cases} \quad (1h)$$

$$P_{i,s}^t = P_{i,Load}^t - P_{i,GT}^t - P_{i,DG}^t - P_{i,dis}^t + P_{i,cha}^t \quad (1l)$$

$$P_{i,w}^{t,\min} \cdot D_{out}^t \leq P_{i,w}^t \leq P_{i,w}^{t,\max} \cdot D_{in}^t \quad (1m)$$

$$D_{out}^t + D_{in}^t \leq 1 \quad (1n)$$

$$\varepsilon_i \cdot \sum_{i \in \mathcal{I}} P_{i,g}^t = P_{i,c}^t \quad (1o)$$

$$\varepsilon \cdot \sum_{i \in \mathcal{I}, t \in T} (P_{i,GT}^t + P_{i,DG}^t) \leq P_c^{\max} \quad (1p)$$

$$P_{i,g}^t = P_{i,GT}^t + P_{i,DG}^t + P_{i,dis}^t - P_{i,cha}^t \quad (1q),$$

where all $t \in T = \{1, \ldots, 24\}$; constraints (1b)-(1d) are the DC power flow model; constraint (1e) are the power limits of transmission lines; (1f) and (1g) denote the output constraint for gas turbines and distributed energy sources; Eq. (1h) represents the energy storage constraint; (1l) is used to calculate the shared electricity; (1m) and (1n) indicates the power limit of interaction with the upstream grid; the carbon emissions are determined by the product of the carbon emission factor and the electrical energy production as given in (1o); (1p) represents the bound on the allowed carbon emissions in the region during 24 hours, and (1q) represents the internal electricity production of VFP $i$.

The RTO determines the shared electricity price and carbon price based on the electricity consumption and carbon emissions of all VFPs.

$$\min_{\lambda_{i,s}, \lambda_{i,c}} \sum_{i \in \mathcal{I}} (\lambda_{i,s}^t)^2 + (\lambda_{i,c}^t)^2 + (\lambda_{i,s}^t - \lambda_{i,s}^{t,k}) + (\lambda_{i,c}^t - \lambda_{i,c}^{t,k}) \quad (2a)$$

$$\text{s.t.} \ \lambda_{i,s}^t = \alpha_i - \beta_i \cdot \sum_{i \in I} P_{i,s}^t \quad (2b)$$

$$\frac{\alpha_i}{\beta_i} - \frac{1}{\beta_i} \lambda_{i,s}^t - \sum_{i \in \mathcal{I}} P_{i,w}^t = 0 \quad (2c)$$

$$-\lambda_{i,s}^t \cdot P_{i,s}^t + \lambda_{i,c}^t \cdot \varepsilon_i = 0 \quad (2e),$$

where objective function (2a) represents the minimization of the price difference between different VFPs; Eq. (2b) is the inverse demand function for electricity; constraint (2c) guarantees the supply and demand equilibrium of the system; (2e) is the first-order condition for carbon prices clearing.

To clear the shared energy market, constraint (2c) ensures system energy balance, which is a common condition. Thus, the market-clearing problem at the transmission can be formulated as GNG and every $VFP_i$ can get its electricity-carbon budget $\Gamma_i$ by investigating the equilibrium of GNG.

### D. VFP sub-model at the Distribution Level

Each VFP allocates its budget within the federation based on transactions at the transmission network level to maximize its social welfare.

*1) Virtual federated consumers:* When the federation is composed of consumer nodes, the optimization problem within the federation can be stated as,

$$\max_{e_n} \sum_{n \in N} \upsilon_n e_n - \frac{\omega_n}{2} e_n^2 - \kappa_n e_n \quad (3a)$$

$$\text{s.t.} \ e_n \in E_n \quad (3b)$$

$$\sum_{n \in N} \kappa_n e_n \leq \Gamma \quad (3c).$$

Here, the objective function is to maximize the social welfare of the coalition consumers, where the first two terms are the consumers' quadratic utility functions, and the last term represents the cost of energy purchase. Constraint (3b) represents the physical constraint satisfied by the user's energy demand; constraint (3c) is that the energy consumption within the VFP meets the budget.

*2) Virtual federated producers:* The optimization problem for a VFP formed spontaneously by all producers is as follows:

$$\max_{e_m} \sum_{m \in M} \kappa_m e_m - a_m e_m^2 - b_m e_m - c_{co2}(e_m) \qquad (4a)$$

$$\text{s.t. } e_m \in E_m \qquad (4b)$$

$$\sum_{n \in N} \kappa_m e_m \leq \Gamma \qquad (4c),$$

which (4a) denotes maximizing the energy supply benefits of the producer coalition, the first term denotes the energy supply benefits, and the last two terms denote the energy supply costs. Constraint (4b) represents the physical constraint satisfied by the producer's energy supply; energy supply benefits need to be kept within budget is given in (4c).

*3) Virtual federated prosumers:* When there are consumers and suppliers in the VFP, the problem of maximizing overall social welfare is expressed as follows:

$$\max_{e_\psi} \sum_{\psi \in N \cup M} \upsilon_n e_n - \frac{\omega_n}{2} e_n^2 - \kappa_\psi e_\psi - c_m(e_m) \qquad (5a)$$

$$\text{s.t. } e_\psi \in E_\psi \qquad (5b)$$

$$e_\psi = e_n - e_m \qquad (5c)$$

$$\sum_{\psi \in N \cup M} \kappa_\psi e_\psi \leq \Gamma \qquad (5d),$$

which $e_\psi$ is the surplus (deficit) of electrical energy within the VFP. Constraint (5b) represents the physical constraint; constraint (5d) is that the energy revenue within the VFP meets the budget.

### E. Budget Allocation Model at Distribution Level

At the distribution network level, the VFP acts as a virtual computing hub platform to adjust budget allocations to achieve optimal global social welfare by issuing weights and benchmark prices based on the electricity information provided by its internal nodes. This section proposes a weight-based budget allocation method that applies to different types of VFPs. First, $VFP_i$ advertises the initial allocation weights $p_n^i/q_m^i$ and benchmark unit prices to the player $n/m$. Then, players send demand/supply $e_n^i / e_m^i$ to $VFP^i$, auctioneer updates the allocation weights according to the demand/supply and republishes the allocation weights. Finally, it adapts to the equilibrium through multiple iterations and an updated total energy surplus $E_i^{k+1} := \sum_{\mathcal{N}_i/\mathcal{M}_i} n/m$ will be available by $VFP_i$.

*1) Consumers:* VFP first releases the allocation weights and benchmark unit electricity prices to each consumer, and then each consumer maximizes its benefits based on all these parameters, and the optimization problem for each consumer is given:

$$\max_{e_n} \upsilon_n e_n - \frac{\omega_n}{2} e_n^2 - p_n \mu_n e_n \qquad (6a)$$

$$\text{s.t. } e_n \in E_n \qquad (6b).$$

When all consumers report their respective demand to VFP, it updates the benchmark price $\mu_n$ based on the demand $e_n$.

$$\mu_n^{k+1} = \mu_n^k - \delta(\sum_{n \in N} p_n \mu_n e_n - \Gamma) \qquad (7).$$

Update allocation weights after the budget are met as follows:

$$p_n^{k+1} = p_n^k + \Delta l * (\Gamma - \sum_{n \in N} p_n \mu_n e_n) / N * \Gamma - p_n^k / \|P\| \qquad (8).$$

The latest weights and benchmark prices are then published for the next optimization until a balanced budget is reached.

*2) Producers:* Each producer gets allocation weights and benchmark prices from VFP to maximize revenue

$$\max_{e_m} q_m \mu_m e_m - a_m e_m^2 - b_m e_m - c_{co2}(e_m) \qquad (9a)$$

$$\text{s.t. } e_m \in E_m \qquad (9b).$$

When all consumers report their supply to VFP, it updates the benchmark price $\mu_m$ based on the supply $e_m$.

$$\mu_m^{k+1} = \mu_m^k - \delta(\sum_{m \in M} q_m \mu_m e_m - \Gamma) \qquad (10).$$

Until the budget is balanced and the allocation weights are updated,

$$q_m^{k+1} = q_m^k + \frac{\Delta l * ((M-2) * \Gamma - \sum_{m \in N} q_m \mu_m e_m)}{(M-1)^2 * \Gamma} - q_m^k / \|Q\| \qquad (11).$$

Stop iterating when the assigned weights no longer change, at which point equilibrium is reached.

*3) Prosumers:* When there are consumers and producers in the VFP, the generated electricity can be treated as negative demand. Thus, the prosumer alliance can be roughly translated into a consumer alliance for the optimal solution.

$$\max_{e_\psi} \sum_{\psi \in N \cup M} \upsilon_n e_n - \frac{\omega_n}{2} e_n^2 - \kappa_\psi e_\psi - c_m(e_m) \qquad (12a)$$

$$\text{s.t. } e_\psi \in E_\psi \qquad (12b)$$

$$e_\psi = e_n - e_m \qquad (12c)$$

When all consumers report their demand to VFP, it updates the benchmark price $\mu_n$ based on the demand $e_n$.

$$\mu_\psi^{k+1} = \mu_\psi^k - \delta(\sum_{\psi \in N \cup M} p_\psi \mu_\psi e_\psi - \Gamma) \qquad (13).$$

After the budget has been balanced, the allocation weights are updated as follows:

$$p_\psi^{k+1} = p_\psi^k + \frac{\Delta l * (\Gamma - \sum_{\psi \in N \cup M} p_\psi \mu_\psi e_\psi)}{(N+M) * \Gamma}) - p_\psi^k / \|P\| \qquad (14)$$

The latest weights and benchmark prices are then published for the next optimization until a balanced budget is reached.

## III. PROPERTIES OF THE GAME

As can be seen in the previous section, the interaction process between VFP and RTO for participating in the energy sharing market at the transmission network layer can be captured by the GNG, and the specific game model is shown in Section A. At the lower distribution network, the stakeholders in VFP compete for budget shares, which is a process that can be characterized by a Nash game regardless of the type of stakeholders involved in the budget allocation, and the specific game model is shown in Section B. The main difference between the GNG and the Nash game is the existence of a common constraint [12].

## A. Game Equilibrium at Transmission Level

The optimization problem of VFP at the transmission network level (1) and the RTO clearing problem (2) together constitute the GNG (2).

1) Players: The set of players $\Upsilon$ consists of VFP $i = 1, \ldots, I$ and trading platform indexed by $i = I + 1$;

2) Strategies: The strategy set of VFP is $\mathcal{P}_i(\lambda_s, \lambda_c) := \{(P_{i,s}, P_{i,c}, P_{i,w}) | (1b) - (1n) \text{ are satisfied.}\}$; the set of TSO is $\mathcal{R}_{I+1}(P_s, P_c, P_w) := \{(\lambda_s, \lambda_c) | (2b) - (2e) \text{ are satisfied.}\}$.

3) Payoffs: The utility function of VFP is $\mathcal{R}_i(\lambda_s, \lambda_c) := \lambda_s \cdot P_{i,s} + \lambda_c \cdot P_{i,c} + \lambda_w \cdot P_{i,w}$, and $\mathcal{R}_{I+1} = \sum_{i \in I} \lambda_s^2 + \lambda_c^2 + (\lambda_s - \lambda_s^k) + (\lambda_c - \lambda_c^k)$.

Thus, the game at the transmission network level is simply expressed as $\mathcal{G}(\Upsilon, \mathcal{P}, \mathcal{R})$. The specific definitions are as follows.

**Definition 1.** The following statements (15) and (16) hold if and only if $(P_s^*, P_c^*, P_w^*, \lambda_s^*, \lambda_c^*) \in \mathcal{P}$ is an equilibrium solution of the generalized Nash game $\mathcal{G}(\Upsilon, \mathcal{P}, \mathcal{R})$ at the transmission level.

$$\left(P_s^*, P_w^*, P_c^*\right) = \underset{P_s, P_c, P_w}{\arg \min} \left\{ \aleph_i(\lambda_s^*, \lambda_c^*), \forall (P_s, P_c, P_w) \in \mathcal{P}_i^* \right\} \quad (15)$$

$$\left(\lambda_s^*, \lambda_c^*\right) = \underset{\lambda_s, \lambda_c}{\arg \min} \left\{ \aleph_{\mathcal{I}+1}, \forall (\lambda_s, \lambda_c) \in \mathcal{P}_{\mathcal{I}+1}^* \right\} \quad (16)$$

In contrast to the traditional Nash game, each player's strategy depends on the rival players' strategies, and its objective function is also related to the variables of the other participants. Hence, the existence and uniqueness conditions for GNG are more stringent.

**Theorem 1.** The solution set $(P_s^*, P_c^*, P_w^*, \lambda_s^*, \lambda_c^*)$ is the equilibrium of the game when and only when $(P_s^*, P_c^*, P_w^*)$ satisfies as the unique solution of (1) and $(\lambda_s^*, \lambda_c^*)$ equal the corresponding dual variable.

**Proof of the generalized Nash game:** The proof can be found in [29].

## B. Game Equilibrium at Distribution Level

The VFP internal budget allocation problem at the distribution level can be regarded as a resource game between rational individuals and is characterized by the Nash game.

1) Virtual federated consumers: From (6), the energy consumption of the consumer $n$ affect the utility function of rival players' energy consumption and the game $\Theta(\mathcal{N}, \Omega_N, \mathcal{F})$ in virtual federated consumers is as follows:

(1) Players: Consumers in set $\mathcal{N}$.

(2) Strategies: According to the weight $p$ received, each consumer determines its energy consumption to maximize energy value. The strategy set of the consumer $n$ is

$$\Omega_N = \left\{ e_n \cdot e_{-n} \middle| e \in E_N \right\}. \quad (17)$$

(3) Payoffs: $\mathcal{F}(e_n; e_{-n})$ is the payoff of consumer $n$, as in (6).

**Definition 2.** When each consumer satisfies $\mathcal{F}(e_n; e_{-n}) \geq \mathcal{F}(e_n'; e_{-n}), \forall e_n' \in \Omega_n$, demand $e_n \in \Omega_n$ is the NE of the game $\Theta(N, \Omega_N, \mathcal{F})$.

At NE points, no entities can unilaterally adjust their bids to obtain higher utility. The nature of NE is analyzed next, as well as the conditions for the existence and uniqueness of NE in the game $\Theta(\mathcal{N}, \Omega_N, \mathcal{F})$.

**Lemma 1.** For the game $\Theta(\mathcal{N}, \Omega_N, \mathcal{F})$, if demand $e_n$ is the NE, let the benchmark unit price be $\mu_n = \Gamma_N / \sum_N e_n$, and the equation $\sum_N \mu_n \cdot e_n \cdot p_n = \Gamma$ holds. Conversely, if there exists $\mu_n > 0$ and the equation $\sum_N \mu_n \cdot e_n \cdot p_n = \Gamma$ holds, and $e_n$ is a NE of the game $\Theta(\mathcal{N}, \Omega_N, \mathcal{F})$.

Based on Lemma 1, the following theorem proves the existence of a unique NE for games $\Theta(\mathcal{N}, \Omega_N, \mathcal{F})$.

**Theorem 2.** For an arbitrary game $\Theta(\mathcal{N}, \Omega_N, \mathcal{F})$, there exists a unique NE $e_n^* = (e_1^*, \ldots, e_N^*) \in \Omega_N$. The budget allocation $\gamma^* = (\gamma_1^*, \ldots, \gamma_N^*)$ at NE satisfies $\gamma_n^* = \mu_n^* \cdot e_n \cdot p_n, \forall n \in \mathcal{N}$, where $\mu_n^* = \Gamma_N / \sum_N e_n^*$ is the corresponding benchmark unit price.

**Proof of the Nash game:** The proof can be found in [29].

2) Virtual federated producers: From (4) it follows that the supply of producer $m$ is related to the supply strategies of other rivals and the game $\Theta(\mathcal{M}, \Omega_M, Z)$ in virtual federated producers is given:

(1) Players: Producers in set $\mathcal{M}$.

(2) Strategies: The energy supply strategy for each customer is influenced by the weights and benchmark prices. The strategy set of the producer $m$ is

$$\Omega_M = \left\{ e_m \cdot e_{-m} \middle| e \in E_M \right\}. \quad (17)$$

(3) Payoffs: $Z(e_m; e_{-m})$ is the payoff of producer $m$, as in (9).

**Definition 3.** When each producer satisfies $Z(e_m; e_{-m}) \geq Z(e_m'; e_{-m}), \forall e_m' \geq 0$, supply $e \in \Omega_M$ is the NE of the game $\Theta(\mathcal{M}, \Omega_M, Z)$.

**Lemma 2.** For the game $\Theta(\mathcal{M}, \Omega_M, Z)$, if supply $e_m$ is the NE, let the benchmark unit price be $\mu_m = \Gamma_M / \sum_M e_m$, and the equation $\sum_M \mu_m \cdot e_m \cdot q_m = \Gamma_M$ holds. Conversely, if there exists $\mu_m > 0$ and the equation $\sum_M \mu_m \cdot e_m \cdot q_m = \Gamma_M$ holds, and $e_m$ is a NE of the game $\Theta(\mathcal{M}, \Omega_M, Z)$.

Based on Lemma 2, the following theorem proves the existence of a unique NE for games $\Theta(\mathcal{M}, \Omega_M, Z)$.

**Theorem 3.** For any game $\Theta(\mathcal{M}, \Omega_M, Z)$, there exists a unique NE $e_m^* = (e_1^*, \ldots, e_M^*) \in \Omega_M$. Moreover, the budget allocation vector $z^* = (z_1^*, \ldots, z_M^*)$ satisfies $z_m^* = \mu_m^* \cdot e_m \cdot q_m, \forall m \in \mathcal{M}$, where $\mu_m^* = \Gamma_M / \sum_M e_m^*$ is the corresponding benchmark unit price.

**Proof of the Nash game:** The proof can be found in [29].

3) Virtual federated prosumers: Corresponding to the model in the previous section, this can be translated into virtual federated consumers to prove.

## IV. IMPLEMENTATION OF MECHANISMS

The joint electricity and carbon market clearing problem are decomposed into the GNG problem at the transmission level and the allocation problem at the distribution level. We solve the two sub-problem alternately to research an equilibrium.

## A. First-order Response Algorithm at Transmission Level

At the transmission level, the electricity price and carbon price of each VFP are obtained with fixed energy $E_i$. A first-order response algorithm (FRA) is proposed to exploit the GNG problem due to the coupling constraints between VFPs. This algorithm alternately optimizes. (1) and (2), where the objective function of (2) reduces the price difference during the iteration by squaring the term.

---

**Algorithm 1**: First-order Response Algorithm

1: Initialization: $k=0$, $e^0$, $P_s^0$, $P_{w,in}^0$, $P_{w,out}^0$, $\varepsilon$;
2: **Repeat**
3:     Platform:
4:     Solve the problem (2) using $P_s^0$, $P_{w,in}^0$, $P_{w,out}^0$ to drive $\lambda_s^k$, $\lambda_c^k$.
5:     VFP $i = 1,...,I$:
6:     Solve the problem (1) using $\lambda_s^k$, $\lambda_c^k$, $e^0$ and drive $P_s^k$, $P_{w,in}^k$, $P_{w,out}^k$.
7:     $k = k + 1$
8: **Until** $|P_s^k - P_s^{k-1}| < \varepsilon$
9: Output GNG equilibrium.

---

Algorithm 1 is capable of converging to equilibrium after several iterations when the condition $\partial^2 h_i/\partial P_{i,g}^2 \geq \beta_i$ is satisfied according to the Cobweb model [30]. From (2b), we derive $\sum P_{i,s} = \alpha_i - 1/\beta_i \cdot \lambda_s^k$, which in economic theory represents the consumer demand curve function is represented by the blue curve in Fig 2, where $\beta_i$ denotes the slope of the demand curve. The cost function $h_i$ in (1a) represents the producer's supply curve function, as indicated by the yellow curve in Fig. 2, where $\partial^2 h_i/\partial P_{i,g}^2$ is the slope of the supply curve. The market converges to the equilibrium situation after several iterations as in Fig. 2(a), which needs to satisfy $\partial^2 h_i/\partial P_{i,g}^2 > \beta_i$; Fig. 2(b) shows the equilibrium fluctuating when $\partial^2 h_i/\partial P_{i,g}^2 = \beta_i$, and finally, Fig. 2(c) shows the non-convergence situation, when $\partial^2 h_i/\partial P_{i,g}^2 < \beta_i$.

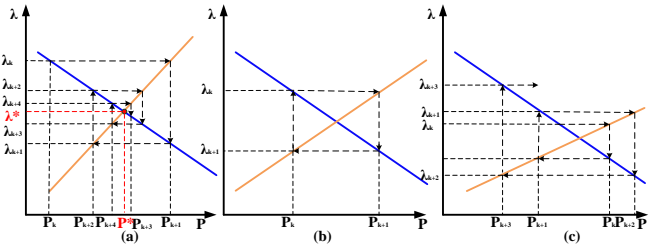

Fig. 2. Schematic diagram of convergence conditions

### B. Distributed Feedback Algorithm at Distribution Level

At the distribution level, the electricity-carbon budget $\Gamma_i$ is treated as an auction lot, allocated by the trusted auctioneer $VFP_i$ through the bidding of the stakeholders $w \in N_i \cup M_i$. Stakeholder $w$ as a rational individual prefers to maximize its profit. However, the individual selfish behavior in the auction will inevitably lead to the loss of the overall social welfare. Thus, a distributed feedback allocation algorithm (DFAA) is devised to reach NE, which combines players' interests with overall social welfare, while the energy surplus $E_i$ of $VFP_i$ is updated.

---

**Algorithm 2**: Distributed Feedback Algorithm

1: Initialization: $k=0$, $\mu^0$, $p^0$, $q^0$, $\varepsilon$;
2: **Repeat**
    VFP advertises $(\mu_N^0, p^0)$ to the consumers;
    // advertises $(\mu_M^0, q^0)$ to the producers;
3:     Consumers solve (6) using $\mu_N^0$, $p^0$ to drive $e_n^k$;
    //producers solve (9) using $\mu_M^0$, $q^0$ to drive $e_m^k$;
4:     Players submit $e^k$ to VFP

---

5:     VFP update $\mu_N^{k+1}$ according (10) to consumers
6:     // update $\mu_M^{k+1}$ according (13) to producers
    $k = k + 1$
7: **Until** $|z^{k+1} - \Gamma_M^{k+1}| < \varepsilon$ and $|\gamma^{k-1} - \Gamma_N^{k+1}| < \varepsilon$
8:     VFP update $p^{k+1}$ according (11)
    // update $q^{k+1}$ according to (14)
    if $\|q^{k+1} - q^k\|_1 < \varepsilon$ and $\|p^{k+1} - p^k\|_1 < \varepsilon$ output
9: Output Nash equilibrium.

---

In general, it is difficult for agents to effectively reach a NE without the assistance of a market algorithm. For this reason, we design a distributed feedback allocation algorithm that enables agents and VFPs to jointly operate the market and reach a NE in a short time. Based on the analysis of the agent's objective function and the nature of NE, the benchmark unit price and weights of energy are chosen as the market signals to regulate the behavior of the agent. When the market signal makes the sum of the agents' budget responses $z^k / \gamma^k$ equal to the total budget $\Gamma^k$, we can determine that the market has converged to the equilibrium point.

### C. Solution of Hierarchical Model

The energy demand and supply within the VFP $E_i$ obtained from the distribution level are then transferred to the transmission level for the next iteration. Global optimum is achieved when the budget of each VFP at each period is close enough between two iterations. The solution producer is listed as:

---

**Algorithm 3**: Adaptive Distributed Alternate Algorithm

1: Initialization: $\tau = 0$, $e^0$, $P_s^0$, $P_{w,in}^0$, $P_{w,out}^0$, $\varepsilon$
2: **Repeat**
3:     Solve the transmission-level GNG problem using **Algorithm 1**, and drive $\Gamma^\tau$.
4:     Each VFP publishes budget to internal nodes
5:     Solve the distribution-level budget allocation problem using **Algorithm 2**, and drive $E^\tau$.
6:     $\tau = \tau + 1$
7: **Until** $|\Gamma_i^\tau - \Gamma_i^{\tau-1}| < \varepsilon$
8: Output equilibrium.

---

When the assumption d) is satisfied, supplied power satisfied $P_{i,g}^{min} + P_{i,w}^{min} < \sum_{i \in} P_{i,load} < P_{i,g}^{max} + P_{i,w}^{max}$, which makes a unique correspondence between the upper and lower layers, thus ensuring the convergence of Algorithm 3.

## V. PERFORMANCE EVALUATION

The designed energy sharing mechanism with budget allocation was tested in a real 110 kV network in Jiangsu province, China, with 2 substations and 16 nodes divided into 5 virtual areas that play the role of VFPs. VFP1 and VFP3 consist entirely of consumer nodes, VFP2 consists of two nodes containing energy storage, two nodes of gas turbine units spontaneously form VFP4, and VFP5 is formed by three consumer nodes and the PV producer node. The structure of the test node system is shown in Fig. 3.

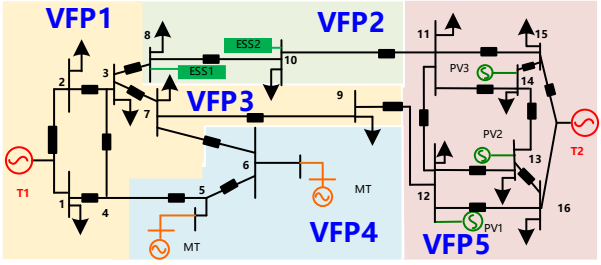

Fig. 3. The structure of the test node system

The wholesale market price, node load data and DER output are all taken from the SGCC for the 12th of July 2020. The specific parameters of the market players can be found in [23]. The simulation will be accomplished on a 64-bit laptop with 4GB-RAM and i5-4210M CPU-2.6GHz using CLPEX and YALMIP.

The analysis of cases mainly consists of the following parts: energy sharing market at the transmission level, budget allocation at the distribution level, analysis of the comparison of mechanisms, and algorithm efficiency.

### A. Analysis of Energy Sharing Market

At the transmission grid level, VFP participates in the wholesale market where the shared energy market and the upper grid interact, and Fig. 4 shows the results of the shared energy market clearing.

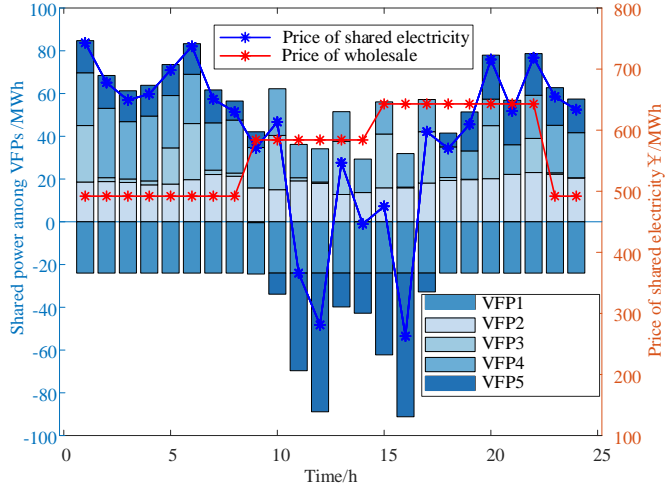

Fig. 4. The results of the shared energy market clearing

When the shared energy is positive, it indicates that the VFP purchases energy from the external world. Conversely, when the shared energy is negative, it indicates that the VFP sells electricity to the external world. Between 8:00 and 18:00, the supply in the shared energy market is greater than the demand due to the higher PV output in VFP5, and the price of shared electricity is lower than the wholesale market price to incentivize market participants to consume electricity internally. As for 20:00 to 7:00 am, the demand in the shared energy market is much higher than the supply, and the price of shared electricity is higher, which causes VFP to purchase much power from the upstream grid to meet the load demand.

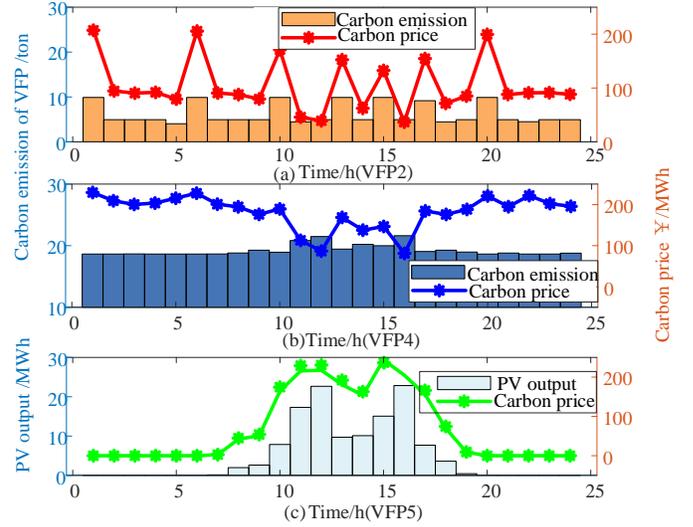

Fig. 5. The results of the carbon market clearing

Three stakeholders are involved in carbon trading at the transmission grid level, including VFP2 with energy storage, VFP4 with combustion turbines and VFP5 with PV. PV in VFP5, as a producer of electrical energy, provides clean energy without generating carbon emissions. Hence, it sells its carbon credits to conventional energy producers, such as VFP2 and VFP4. As shown in Fig. 5, the higher the PV output, the higher its sold carbon emissions price.

The gas turbine in VFP4 has the highest carbon emission factor as the most polluting equipment in the system. To meet the carbon credits of VFP4, carbon credits need to be purchased from the remaining VFPs. As a result, VFP4 has the highest carbon price among the three. The increase in clean energy-PV output between 10:00 and 16:00 leads to an increase in carbon credits circulating in the market and a relatively low price for purchasing carbon credits in VFP4, which makes VFP4 increase its carbon emissions during this period.

Depending on its internal energy storage operations, VFP2 purchases carbon credits from other participants when carbon emissions exceed the allowances and resells the surplus to participants who need it when there is a surplus of carbon emissions. Due to the limited carbon credits, when VFP2 has small carbon emissions, it will sell carbon credits at a lower price to obtain market share; in contrast, when its carbon emissions are high, it sets a higher carbon price to protect its emissions.

### B. Analysis of Budget Allocation Results

The power output of all devices in the system is shown in Fig. 6 after several iterations of interaction between the upper transmission network and the lower distribution network to reach equilibrium. When $P_w$ is less than 0, it indicates that the system sells electricity to the upstream grid. A negative output of ESS indicates charging, and a positive one indicates discharging.

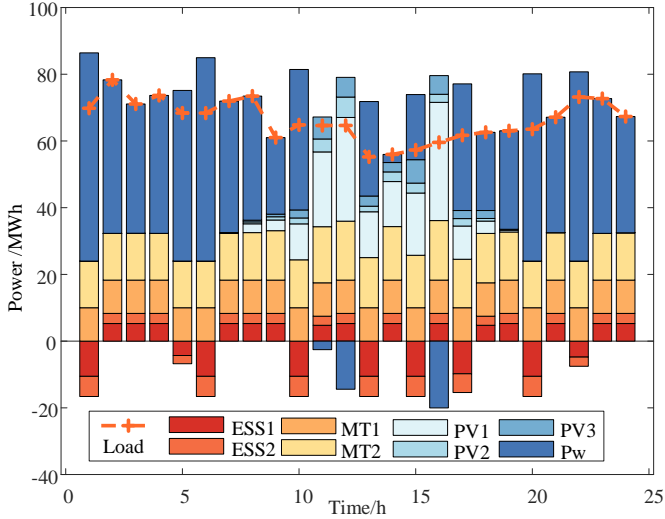

Fig. 6. The power output of all the equipment in the system

Table I shows the results of the internal budget allocation for the VFP containing the load. Both VFP1 and VFP3 are federations of load nodes, and the higher the energy use value of a node, the more budget it is allocated. VFP5, as a coalition of prosumers, does not follow the above variation law of pure load nodes for nodes 12,13,14 containing PV. The budget allocated at this point indicates profit. Although node 12 has the highest value of energy use, its PV output is not as much power as node 13, so its profit allocation is lower than that of node 13. Node 13 has the lowest load energy use utility, but its PV output is the highest, getting the highest share of the profit.

TABLE I
RESULTS FOR BUDGET ALLOCATION

| VFP | Nodes | $\Gamma$(¥) | Value of energy use $\upsilon$ (¥) | |
|---|---|---|---|---|
| VFP1 | 1 | 1292.5 | 1,808.6 | |
| | 2 | 4.83 | 2.8 | |
| | 3 | 1271.9 | 1,716.1 | 3,527.6 |
| | 4 | 0 | 0 | |
| VFP3 | 7 | 1338.3 | 888.3 | |
| | 9 | 1449.28 | 1,780.4 | 2,668.7 |
| VFP5 | 11 | 145.4 | 16.0 | |
| | **12** | **961.2** | **1,049.9** | |
| | **13** | **1289.7** | **264.25** | |
| | **14** | **741.3** | **40.0** | 1647.4 |
| | 15 | 1126.1 | 253.1 | |
| | 16 | 154.3 | 24.02 | |

## C. Analysis of the Comparison of Mechanisms

In this section, the proposed energy sharing mechanism with a budget allocation is contrasted with energy sharing only without budget allocation and energy sharing with budget allocation, respectively.

TABLE II
COMPARATIVE OF CONSUMERS UNDER DIFFERENT MECHANISMS

| | Value of energy (¥) | Flexible demand (MWh) |
|---|---|---|
| **Energy sharing with budget allocation** | **7,843.7** | **110.64** |
| Energy sharing without budget allocation | 7022.5 | 184.4 |
| Without energy sharing and budget allocation | 6,183.3 | 184.4 |

As can be seen from Table II, although the total amount of electricity involved in demand response in the method proposed in this paper is not as much as in the other two cases, the consumer's energy use value is the greatest under this method. Our proposed method proves to be very effective in improving energy efficiency.

TABLE III
COMPARATIVE OF PRODUCERS UNDER DIFFERENT MECHANISMS

| | | Cost (¥) | Carbon Emission (t) |
|---|---|---|---|
| **Energy sharing with budget allocation** | GT | **3,882.36** | **462.01** |
| | ESS | **37.73** | **3.85** |
| Energy sharing without budget allocation | GT | 6,216.9 | 733.230451 |
| | ESS | 10.031 | 1.02 |
| Without energy sharing and budget allocation | GT | 6,637.2 | 1701.16 |
| | ESS | 26.94 | 2.75 |

Table III compares the costs and carbon emissions of the system electrical energy production equipment for the three scenarios. With budget allocation, there is a significant reduction in the cost of GT and a significant reduction in GT carbon emissions. While the cost of ESS without budget allocation is relatively lower, the carbon emission factor using ESS is lower than GT. Therefore, the overall view is that using budget allocation can better reduce the carbon emissions of the system.

## D. Analysis of Computational Efficiency

The number of iterations and solution time of the proposed distributed solution algorithm are shown in Table IV. The outer iteration corresponds to Algorithm 3 to realize the transmission and distribution network coupling allocation, the transmission level iteration corresponds to Algorithm 1, and the distribution network level budget allocation corresponds to Algorithm 2. When there are more participants, the number of iterations of Algorithm 2 increases, and the solution time takes longer.

TABLE IV
COMPUTATIONAL PERFORMANCE IN ITERATION PROCESS

| Computational performance | | Iteration | Solution time (s) |
|---|---|---|---|
| Outer iteration | | 3 | 263.01 |
| Transmission-level | | 3 | 63.59 |
| Distribution-level (Budget allocation) | VFP1 | 19 | 109.17 |
| | VFP2 | 34 | 98.48 |
| | VFP3 | 14 | 6.46 |
| | VFP4 | 16 | 7.85 |
| | VFP5 | 58 | 199.42 |

Although the overall computation time is 263.01 seconds, the result is for the previous day's 24 cases, not a time section, so the solution time is acceptable. The case study shows that the proposed algorithm has good convergence.

## VI. CONCLUSION

This paper proposes a VFP-based hierarchical market framework that combines transmission and distribution networks within a region, with simultaneous clearing electricity and carbon market. At the transmission level, VFPs participate as stakeholders in energy sharing markets, wholesale markets and carbon markets, whereas at the distribution level, act as a selfless auctioneer to allocate the energy budget. The interaction among VFPs at the transmission level is explained as GNG, for which a first-order optimal response algorithm is proposed to clear the energy sharing market. The Nash game is employed to characterize the budget allocation problem at the distribution network level, and a distributed feedback allocation algorithm is developed to obtain equilibrium. Case studies show that the energy sharing mechanism with budget allocation provides a good incentive for DERs to participate in the market to improve consumer energy efficiency and reduce carbon emissions. Furthermore, the existence and uniqueness of the game were also proved in detail. The effectiveness and convergence of the algorithm are also well proven.